\documentclass[10pt]{amsart}
\usepackage{amssymb,amsmath,amsthm,amsfonts,amsopn,url}
\usepackage{amscd,amssymb,amsopn,amsmath,amsthm,graphics,amsfonts,enumerate,verbatim,calc}
\usepackage[all]{xy}
\theoremstyle{plain}
\newtheorem{thm}{Theorem}[section]
\newtheorem*{mt*}{Theorem (Huneke-Lyubeznik)}
\newtheorem*{cj*}{Conjecture}
\newtheorem*{nt*}{Notations}

\newtheorem{lemma}[thm]{Lemma}
\newtheorem{cor}{Corollary}

\theoremstyle{definition}

\newcommand{\ideal}[1]{\mathfrak{#1}}

\newcommand{\n}{\ideal{n}}

\newcommand{\func}[1]{\mathrm{#1} \,}

\newcommand{\height}{\func{ht}}

\newcommand{\Ass}{\func{Ass}}

\newcommand{\ZZ}{{\mathbb Z}}

\title[]{Finiteness of local cohomology modules over arbitrary regular local rings}

\author[]{Rajsekhar Bhattacharyya}

\address{Dinabandhu Andrews College, Garia, Kolkata 700084, India}

\email{rbhattacharyya@gmail.com}

\thanks{}

\keywords{Local Cohomology}

\subjclass[2010]{13D45}
\begin{document}

\begin{abstract}
In this paper, we consider arbitrary regular local ring $S$ and an ideal $J\subset S$ whose minimal primes are cohomologically complete intersection. In Section 2, using the Mayer-Vietoris spectral sequence, we construct examples of local cohomology modules $H^i_J (S)$, whose set of associated primes is finite, where minimal primes satisfy the `Dropping Condition'. As a special case, we show that when minimal primes of $J$ are part of a specific regular system of parameters satisfy the `Dropping Condition', the set of associated primes of $H^i_J (S)$ is finite. 

In Section 3, we prove the main result without the `Dropping Condition' and as a special case, we show that when minimal primes of $J$ are part of a specific regular system of parameters, the set of associated primes of $H^i_J (S)$ is finite. 

These examples of local cohomology modules include the cases whose associated primes can contain prime number $p$, but still the set is finite.
\end{abstract}

\maketitle
\section{introduction}

Let $R$ be a Noetherian ring and $M$ be a module over it. For an ideal $I\subset R$ and for some integer $i\geq 0$, consider local cohomology module $H^i_I (M)$ with support in the ideal $I$. In the fourth problem of \cite{Hu}, it is asked that whether local cohomology modules of Noetherian rings have finitely many associated prime ideals. There are examples given in \cite{Si}, \cite{Ka}, and \cite{SS}, which show that the set of associated primes of $i$-th local cohomology module $H^i_I (R)$ of a Noetherian ring $R$ can be infinite. However, there are several important cases where we have the finiteness of the set of associated primes of local cohomology modules and we list them below: (1) Regular rings of prime characteristic \cite{HS}, (2) regular local and affine rings of characteristic zero \cite{Ly1}, (3) unramified regular local rings of mixed characteristic \cite{Ly2} and (4) smooth algebra over $\ZZ$ \cite{BBLSZ}. On the other direction, the result is also true for every low dimensional regular local ring upto dimension five \cite{Boix}. These results support Lyubeznik conjecture, (\cite{Ly1}, Remark 3.7): 

\textbf{Conjecture:} Let $R$ be a regular ring and $I\subset R$ be its ideal, then for each $i\geq 0$, $i$-th local cohomology module $H^i_I (R)$ has finitely many associated prime ideals.

Recently, in \cite{Ma} it is shown that in ramified regular local rings, there exist examples of local cohomology modules with ifinitely many associated primes. To investigate the cause in more detail, in this paper, we proceed as follows, and it is to be noted that the results here, are true for any regular local ring. 

At first, in Section 2 (The Main Result I), we consider arbitrary regular local ring $S$ and an ideal $J\subset S$ whose minimal primes are cohomologically complete intersection \cite{HeSc} and satisfying the `Dropping Condition' (see definition inside). Using the Mayer-Vietoris spectral sequence, we construct examples of local cohomology modules $H^i_J (S)$, whose set of associated primes is finite. For detail, see, Theorem 2.3. As a special case, in Corollary 1, we show that when minimal primes of $J$ are part of a specific regular system of parameters satisfy the `Dropping Condition', the set of associated primes of $H^i_J (S)$ is finite. This result is consistent with \cite{Ma}, since the ideal there does not satisfy the `Dropping Condition', for detail see Remark 2.

In Section 3 (The Main Result II), we prove Theorem 2.3 and Corollary 1, in more general sence, without the use of `Dropping Condition'. See Theorem 3.1 and Corollary 2. Corollary 2 tells that for an ideal $J\subset S$ if its minimal primes are part of one regular system of parameters, then the set of associated primes of $H^i_J (S)$, is finite. This is not consistent with the result obtained in \cite{Ma}

Finally, in \cite{Nu3}, it is shown that if the set of associated primes of a local cohomology module does not contain the prime number $p$, then the set of associated primes is finite. These examples of local cohomology modules (in this paper) include the cases whose associated primes can contain prime number $p$, but still the set is finite, see Remark 5.

\section{the main result I}






In proving our theorem, we need the Mayer-Vietoris spectral sequence, which arises from te Roos complex \cite{Roo} where in \cite{AGZ}, it is stated and proved only for defining ideal of an arrangement of linear subvarieties. In \cite{Ly4} it is extened for arbitrary commutative ring and module where \cite{ABZ}, the spectral sequnce has been constructed and properties are discussed very elaborately for more general categorical objects. Although the Mayer-Vietoris spectral sequences as observed in both \cite{AGZ} and \cite{Ly4}, differ in $E_0$ and $E_1$ page, but the coincide on $E_2$ page (see 4.2 Example \cite{ABZ} for detail).  

In this paper, we need the following theorem from \cite{Ly4} regarding the Mayer-Vietoris spectral sequence. 

\begin{thm}[Theorem on MVSS]\label{MV}
Let $A$ be an arbitrary commutative ring, let $I_1,\dots,I_n\subset A$ be ideals and
let $M$ be an $A$-module. There exists a spectral sequence
$$E_1^{-p,q}=\bigoplus_{i_0<\dots<i_p}H^q_{I_{i_0}+\dots+I_{i_p}}(M)
\Longrightarrow H^{q-p}_{I_1\cap\dots\cap I_n}(M).$$
\end{thm}

In our theorem we need the condition so that the above spectral sequence will be degenarate on $E_2$-page. In \cite{AGZ} in Theorem 1.2 it is shown the degenaracy in $E_2$-page for a defining ideal of an arrangement of linear subvarieties, when the underlying regular ring contains field. In \cite{ABZ} (see Lemma 4.15), degenaracy condition is analysed in a general setting and ideal which is cohomologically complete intersection (i.e. local cohomology module in the support of the ideal is non zero only at its height) used for the degenaracy. In \cite{Ly4} it has been mentioned that the spectral sequence will be degenarate on $E_2$ page if $A$ and $A/I_i$ are all regular (possibly arbitrary), with $M=A$ and gives reference of Theorem 1.2 of \cite{AGZ}.

To prove the main theorem, we also need minimal primes of the support ideal to be cohomologicaly complete intersection whose sums are also so. For an ideal $I$ in regular local ring $S$, if $S/I$ is regular then $I$ is cohomologicaly complete intersections (although converse is not true, consider $k[[x,y,z]]$ and the ideal $(x^2 +y^3 + z^5)$). 

\begin{lemma}
Let $(A,\o)$ regular local ring.\\ 
(a) An ideal $K$ will be generated by regular system of parameters, or equivalently $A/K$ is regular if and only if $K\cap {\o}^2\subset \o K$.\\
(b) If $I, J$ are ideals generated by part of regular system of parameters, or equivalently $A/I, A/J$ are regular rings and $A/(I+J)$ is also a regular then regular system of parameter of $I+J$ can get obtained by joining part of regular system of parameters of $I$ to that of $J$ and conversely if $I$ and $J$ are formed from a part of regular system of parameters then $I+J$ is generated by a regular system of parameters. 
\end{lemma}

\begin{proof}
(a) The asssertion will be true if and only if the minimal set of generators of $K$ is a regular system of parameters, since a regular system of parameters is a minimal set of generators for the ideal it generates. Now the minimala set of generators forms basis of the vector space $K/\o K$ over $A/\o$. So $K$ will be genetared by regular system of parameters if and only if the dimension of the vector space $(K+{\o}^2)/{\o}^2\cong K/K\cap {\o}^2$ is same to that of $K/\o K$ over $A/\o$. Since there is already a surjective canonical map from $K/\o K$ to $K/K\cap {\o}^2$, thus $K/\o K$ is canonically isomorphic to $K/K\cap {\o}^2$ if and only if $K\cap {\o}^2\subset \o K$. 

(b) Let $I=({u}_1,\ldots,{u}_k)A$ and $J=({v}_1,\ldots,{v}_l)A$ be generated by regular system of parameters and let $\{{u}_{i_1},\ldots,{u}_{i_a}, {v}_{j_1},\ldots,{v}_{j_b}\}$ be a set of minimal generatos of $I+J$. 
From (a) of this lemma it follows that they will form the regular system of parameters and our assertion follows. Converse follows from the fact that joining two part of same regular system of parameters will be a part of regular system of parameters.
\end{proof}

Now we prove the main result of the paper.

\begin{thm}
Let $(S,\n)$ be a 
complete ramified regular local ring of mixed characteristic and $J$ be an ideal of it. For $t=1,\ldots, n$, let $I_t$'s be the ideals  of $S$ and $J=\cap^n_{t=1} I_t$, such that all the ideals $\sum_{t\in F}I_t$ are prime and cohomologically complete intersection for every subsets $F\subset\{1,\dots,n\}$. Moreover we assume that dropping any $I_t$ from the $\sum_{t\in F}I_t$ gives proper containment. Then $\Ass H^i_{J} (S)$ is finite.
\end{thm}

\begin{proof}

The basic idea of the proof is as follows: Since for every subsets $F\subset\{1,\dots,n\}$, $\sum_{t\in F}I_t$ are all prime and cohomologically complete intersection, the Mayer-Vietoris spectral sequence will degenerate on $E_1$-page. This will imply that $H^i_{J} (S)$ has a finite filtration such that every subquotient module of the filtration is isomorphic to a local cohomology module whose set of associated primes is finite. Thus the set of associated primes of $H^i_{J} (S)$ is finite. 

We apply the above theorem on MVSS for the regular local ring $S$ and ideal $J$ with $M=S$. Here we start from $E_1$ page and define $E_1^{-p,q}$ according to \cite{Ly4}, infact previously, it has actually defined in a more general context in \cite{AGZ}. Set $J_{i_0\ldots i_p}=I_{i_0}+\dots+I_{i_p}$ and $J_{i_0\ldots\widehat{i_j}\ldots i_p}=I_{i_0}+\dots+I_{i_{j-1}}+I_{i_{j+1}}+\dots+I_{i_p}$:
$$E_1^{-p,q}=\bigoplus_{i_0<\dots<i_p}H^q_{J_{i_0\ldots i_p}}(S)
\Longrightarrow H^{q-p}_{I_1\cap\dots\cap I_n}(S)= H^{q-p}_{J}(S)$$

where differentials $E_1^{-p,q}\to E_1^{-p+1,q}$ on $E_1$-page are horizontal and they are defined according to \cite{Ly4} and in the following way: Consider their restriction to $H^q_{J_{i_0\ldots i_p}}(S)\subset E_1^{-p,q}$ where they take elements of it to some local cohomology module $H^{q}_{J_{i_0\ldots\widehat{i_j}\ldots,i_p}}(S)\subset E_2^{-p+1,q}$ and this map is induced by natural inclusion of $J_{i_0\ldots\widehat{i_j}\ldots i_p}\subset J_{i_0\ldots i_p}$. 


If in $E_1$-page if $H^q_{J_{i_0\ldots i_p}}(S)$ is a non zero summand in $E_1^{-p,q}$ then $H^l_{J_{i_0\ldots\widehat{i_j}\ldots i_p}}(S)$ will be non zero for $l\leq q-1$. This is because we assume that dropping an index we have proper containment and since all $J_{i_0\ldots i_p}$ are primes, after dropping one index the height will be dropped by atleast one. Thus any nonzero local cohomology in the support of $J_{i_0\ldots\widehat{i_j}\ldots,i_p}$ appears along or below the diagonal of $H^q_{J_{i_0\ldots i_p}}(S)$.

As defined above, in $E_1$-page differentials are all zero. Thus on the $E_2$ page, we get $E_2^{-p,q}= E_1^{-p,q}=\bigoplus_{i_0<\dots<i_p}H^q_{J_{i_0\ldots i_p}}(S)$ along with the differentials $E_2^{-p,q}\to E_2^{-p+2,q-1}$. Consider a differential from $E_2^{-p,q}\to E_2^{-p+2,q-1}$ on the $E_2$-page and consider its restriction to $H^q_{J_{i_0\ldots i_p}}(S)\subset E_2^{-p,q}$ where they take elements of it to $E_2^{-p+2,q-1}=\bigoplus_{j_0<\dots<j_{p-2}}H^{q-1}_{J_{j_0\ldots j_{p-2}}}(S)$. If $H^{q-1}_{J_{j_0\ldots j_{p-2}}}(S)= 0$, then we have the zero map from $H^q_{J_{i_0\ldots i_p}}(S)$. Assume $H^{q-1}_{J_{j_0\ldots j_{p-2}}}(S)\neq 0$, this implies $\height J_{j_0\ldots j_{p-2}}= q-1$. Since $\height J_{i_0\ldots i_p}= q$, $J_{j_0\ldots j_{p-2}}$ can not contain $J_{i_0\ldots i_p}$ then according to Lemma 4.15 of \cite{ABZ}, there exists only the zero map from $H^q_{J_{i_0\ldots i_p}}(S)$ to $H^{q-1}_{J_{j_0\ldots j_{p-2}}}(S)$. So there is only zero differential in $E_2$-page.

In fact for any positive integer $s$, inductively on $E_s$-th page, we have $E_{1}^{-p,q}= E_s^{-p,q}$ and by the same arguement as given to the $E_2$-page, the differential on the $E_s$-page is also zero. Thus $E_{1}^{-p,q}= E_{s+1}^{-p,q}$. Thus $E_{\infty}^{-p,q}= E_2^{-p,q}$. So for fixed $i=q-p$, $H^{i}_{J}(S)= H^{q-p}_{J}(S)= H^{q-p}_{I_1\cap\dots\cap I_n}(S)$ and we will get a filtration whose subquotient modules are isomorphic to $E_{\infty}^{-p,q}$ with varying $q$ and $p$ such that $i=q-p$. Since at each page of the spectral sequence we have only finitely many rows and columns, we have finite filtration and it is then sufficient to show that for each $(-p,q)$-pair, the set of asssociated primes of $E_{\infty}^{-p,q}$ is finite. From \cite{Mar} we know for every ideal $J_{j_0\ldots j_{p-2}}$ of $S$, $\Ass H^i_{J_{j_0\ldots j_{p-2}}}(S)$ is finite. So for each $(-p,q)$-pair $\Ass E_{\infty}^{-p,q}$ is finite. This implies $\Ass H^i_{J}(S)$ is contained in the subset of $\{J_{i_0\ldots i_p}:q-p=i\}$ and it is finite. This concludes the proof. 
\end{proof}

{\bf Definition}\\

In reference to the statement of Theorem 2.3, if we assume that dropping any $I_t$ from the $\sum_{t\in F}I_t$ gives proper containment, we define it as `Dropping Condition' satisfied by the minimal prime ideals.\\

{\bf Remark 1}\\

We can realize the above theorem in the following situation: For a ramified regular local ring $S$, let $x_1,\dots,x_d$ be a regular system of parameters such that $I_t$'s of the theorem are part of this regular system of parameters. From Lemma 2.2 it follows that the necessary and sufficient condition for the sum of any two such $I_t$ to be a part of regular system of parameters is to take all the $I_t$'s from one such regular system of parameters. Moreover we choose the $I_t$'s such that dropping any $I_t$ from the $\sum_{t\in F}I_t$ for every subsets $F\subset\{1,\dots,n\}$, gives proper containment. As for example for $n=3$, $I_1=(x_1,x_2)$, $I_2=(x_1,x_3)S$ and $I_3=(x_1,x_4)S$ taken from the regular sysytem of parameters $x_1,\ldots, x_d$, satisfy this condition. Thus we get the following corollary.

\begin{cor}
Consider a part of regular system of parameters $x_{i_1},\dots,x_{i_k}$ of a ramified regular local ring and consider ideals $I_t$ for $t=1,\ldots, n$ such that they are made out of parameters arbitrarily. Moreover we choose the $I_t$'s such that they satisfy the `Dropping Condition'. Then for the ideal $J=\cap^n_{i=1}I_t$ the set of associated primes of $H^i_{J}(S)$ is finite.
\end{cor}

{\bf Remark 2}\\

In \cite{Ma}, for $R = S/(2 + y^2_1 y^2_3 y_6 + y_1y_2y_3y_4y_5 + y^2_2 y^2_4 y_6)$ where $S = \mathbb{Z}_2[[x_1,\ldots, x_6, y_1,\ldots, y_6]] $ and for the ideal $J=IS +(y_1, y_2)S$ where $I\subset \mathbb{Z}_2[[x_1,\ldots, x_6]]$ is the square free monomial ideal corresponds to six vertex triangulation of $\mathbb{R}\mathbb{P}^2$, it is shown that local cohomology $H^6_J (R)$has infinitely many associated primes. Incidentally the set of minimal prime ideals of $J$ are $\{(\{x_i: i\notin H\}, y_1, y_2)\}$ (where $H$ varies over facets of the triangulation of $\mathbb{R}\mathbb{P}^2$), which part of regular system of parameters but the do not satisfy the `Dropping Condition' as defined above. 

\section{the main result II}

In this section, we prove Theorem 2.3 in more general sence, without the use of `Dropping Condition'. 

\begin{thm}
Let $(S,\n)$ be a 
regular local ring and $J$ be an ideal of it. For $t=1,\ldots, n$, let $I_t$'s be the ideals  of $S$ and $J=\cap^n_{t=1} I_t$, such that all the ideals $\sum_{t\in F}I_t$ are prime and cohomologically complete intersection for every subsets $F\subset\{1,\dots,n\}$. Then $\Ass H^i_{J} (S)$ is finite.
\end{thm}

\begin{proof}

The basic idea of the proof is as follows: Since for every subsets $F\subset\{1,\dots,n\}$, $\sum_{t\in F}I_t$ are all prime and cohomologically complete intersection, the Mayer-Vietoris spectral sequence will degenerate on $E_2$-page. This will imply that $H^i_{J} (S)$ has a finite filtration such that every subquotient module of the filtration is isomorphic to a local cohomology module whose set of associated primes is finite. Thus the set of associated primes of $H^i_{J} (S)$ is finite. 

We apply the theorem on MVSS i.e. Theorem 2.1 above for the regular local ring $S$ and ideal $J$ with $M=S$. Here we start from $E_1$ page and define $E_1^{-p,q}$ according to \cite{Ly4}, infact previously, it has actually defined in a more general context in \cite{AGZ}. Set $J_{i_0\ldots i_p}=I_{i_0}+\dots+I_{i_p}$ and $J_{i_0\ldots\widehat{i_j}\ldots i_p}=I_{i_0}+\dots+I_{i_{j-1}}+I_{i_{j+1}}+\dots+I_{i_p}$ and from Theorem 2.1, we have the following convergence of the spectral sequence:
$$E_1^{-p,q}=\bigoplus_{i_0<\dots<i_p}H^q_{J_{i_0\ldots i_p}}(S)
\Longrightarrow H^{q-p}_{I_1\cap\dots\cap I_n}(S)= H^{q-p}_{J}(S)$$

Here differentials $E_1^{-p,q}\to E_1^{-p+1,q}$ on $E_1$-page are horizontal and they are defined according to \cite{Ly4} and in the following way: Consider their restriction to $H^q_{J_{i_0\ldots i_p}}(S)\subset E_1^{-p,q}$ where they take elements of it to some local cohomology module $H^{q}_{J_{i_0\ldots\widehat{i_j}\ldots,i_p}}(S)\subset E_2^{-p+1,q}$ and this map is induced by natural inclusion of $J_{i_0\ldots\widehat{i_j}\ldots i_p}\subset J_{i_0\ldots i_p}$. 

Let $H^q_Q (R)$ be a nonzero direct summand of $E_1^{-p,q}$ and there is a map from here to $H^q_P (R)$ in $E_1^{-p+1,q}$ induced by the inclusion $P\subset Q$. Clearly $\height Q=q$. If $\height P=q$ then $Q=P$ and the map will be identity, otherwise for $\height P\neq q$, $H^q_P (R)=0$ and the map is a zero map. Now at each summand $H^q_{J_{i_0\ldots i_p}}(S)$ of $E_1^{-p,q}$ we map to $H^q_{J_{i_0\ldots\widehat{i_j}\ldots i_p}}(S)$ where $J_{i_0\ldots\widehat{i_j}\ldots i_p}\subset J_{i_0\ldots i_p}$. If for every $i_j$, $J_{i_0\ldots\widehat{i_j}\ldots i_p}\subset J_{i_0\ldots i_p}$ properly, then $H^q_{J_{i_0\ldots\widehat{i_j}\ldots i_p}}=0$ and the restriction of the differential on $H^q_{J_{i_0\ldots i_p}}(S)$ will be a zero map. Now the differential map which ends on $E_2^{-p,q}$, the image of it in $H^q_{J_{i_0\ldots i_p}}(S)$, is either the whole $H^q_{J_{i_0\ldots i_p}}(S)$ or zero. Since the modules on $E_2$-page are the homology modules of the differentials of $E_1$-page and since homology commutes with the direct sum, in this case either $H^q_{J_{i_0\ldots i_p}}(S)$ appears or will not appear in $E_2^{-p,q}$ as a direct summand. On the other hand, if for any one $i_j$, $J_{i_0\ldots\widehat{i_j}\ldots i_p}= J_{i_0\ldots i_p}$ then the map is identity, which implies kernel of this differential restricted at $H^q_{J_{i_0\ldots i_p}}(S)$ is zero and this particular local cohomology summand module can not appear in $E_2^{-p,q}$. Thus on the $E_2$ page, we get $E_2^{-p,q}\subset E_1^{-p,q}=\bigoplus_{i_0<\dots<i_p}H^q_{J_{i_0\ldots i_p}}(S)$ and the differentials are from $E_2^{-p,q}\to E_2^{-p+2,q-1}$ on the $E_2$-page. Consider its restriction to $H^q_{J_{i_0\ldots i_p}}(S)\subset E_2^{-p,q}$ where they take elements of it to $E_2^{-p+2,q-1}\subset E_1^{-p+2,q-1}= \bigoplus_{j_0<\dots<j_{p-2}}H^{q-1}_{J_{j_0\ldots j_{p-2}}}(S)$. If $H^{q-1}_{J_{j_0\ldots j_{p-2}}}(S)= 0$, then we have the zero map from $H^q_{J_{i_0\ldots i_p}}(S)$. Assume $H^{q-1}_{J_{j_0\ldots j_{p-2}}}(S)\neq 0$, this implies $\height J_{j_0\ldots j_{p-2}}= q-1$. Since $\height J_{i_0\ldots i_p}= q$, $J_{j_0\ldots j_{p-2}}$ can not contain $J_{i_0\ldots i_p}$ then according to Lemma 4.15 of \cite{ABZ}, there exists only the zero map from $H^q_{J_{i_0\ldots i_p}}(S)$ to $H^{q-1}_{J_{j_0\ldots j_{p-2}}}(S)$. So there is only zero differential in $E_2$-page.  
 
In fact for any positive integer $s$, inductively on $E_s$-th page, we have $E_{2}^{-p,q}= E_s^{-p,q}$ and by the same arguement as given to the $E_2$-page, the differential on the $E_s$-page is also zero. Thus $E_{2}^{-p,q}= E_{s+1}^{-p,q}$. Thus $E_{\infty}^{-p,q}= E_2^{-p,q}$. So for fixed $i=q-p$, $H^{i}_{J}(S)= H^{q-p}_{J}(S)= H^{q-p}_{I_1\cap\dots\cap I_n}(S)$ and we will get a filtration whose subquotient modules are isomorphic to $E_{\infty}^{-p,q}$ with varying $q$ and $p$ such that $i=q-p$. Since at each page of the spectral sequence we have only finitely many rows and columns, we have finite filtration and it is then sufficient to show that for each $(-p,q)$-pair, the set of asssociated primes of $E_{\infty}^{-p,q}$ is finite. From \cite{Mar} we know for every ideal $J_{j_0\ldots j_{p-2}}$ of $S$, $\Ass H^i_{J_{j_0\ldots j_{p-2}}}(S)$ is finite. So for each $(-p,q)$-pair $\Ass E_{\infty}^{-p,q}$ is finite. This implies $\Ass H^i_{J}(S)$ is contained in the subset of $\{J_{i_0\ldots i_p}:q-p=i\}$ and it is finite. This concludes the proof. 
\end{proof}

{\bf Remark 3}\\

We can realize the above theorem in the following situation: For a regular local ring $S$, let $x_1,\dots,x_d$ be a regular system of parameters such that $I_t$'s of the theorem are part of this regular system of parameters. From Lemma 2.2 it follows that the necessary and sufficient condition for the sum of any two such $I_t$ to be a part of regular system of parameters is to take all the $I_t$'s from one such regular system of parameters. Thus we get the following corollary.

\begin{cor}
Let $S$ be a regular local ring and $J$ be its an ideal such that its minimal primes $I_t$ for $t=1,\ldots, n$ are part of regular system of parameters of a specific regular system of parameters $x_{1},\dots,x_{d}$ of $S$. Then the set of associated primes of $H^i_{J}(S)$ is finite.
\end{cor}

\begin{proof}
Since the ideal $J=\cap^n_{i=1}I_t$, from Remark 3 and Theorem 3.1 the result follows.
\end{proof}

{\bf Remark 4}\\

From the discussion in Remark 2 it follows that the result of this corollary is not consistent with that obtained in \cite{Ma}.\\

{\bf Remark 5}\\

In \cite{Nu3}, it is shown that if the set of associated primes of a local cohomology module does not contain the prime number $p$, then the set of associated primes is finite. But this is not a necessary condition. We show in the following way:

Consider the ramified regular local ring $S= {\mathbb{Z}}_p[[X_1,\dots, X_d]]/(X^2_1-p)$. Clearly in $S$, $x_1,\dots, x_d$ the images of $X_1,\dots, X_d$, form a regular system of parameters, such that $x^2_1=p$. If we take all the $I_t$'s which contains $x_1$, and satisfy the condition of the Corollary 1, then $J$ contains $x_1$ and hence $x^2_1= p$. From Corollary 1, we find that the set of associated primes of local cohomology modules $H^i_J (S)$ is finite and its associated primes contain $p$.

\section{acknowledgement}

I would like to thank Professor Alberto F. Boix for pointing out errors in the earlier manuscript and providing the counter example that every cohomologically complete intersection ideal is not such an ideal that its quotient ring is regular local.

\end{document}